\documentclass[12pt]{amsart}
\usepackage{graphicx}
\usepackage{amssymb}
\usepackage{amsmath}
\usepackage{amsthm}
\usepackage{amscd}
\newtheorem{thm}{Theorem}[section]
\newtheorem{cor}[thm]{Corollary}
\newtheorem{lem}[thm]{Lemma}

\theoremstyle{definition}
\newtheorem{defn}[thm]{Definition}
\newtheorem{rem}[thm]{Remark}
\numberwithin{equation}{section}

\newcommand{\D}{\Delta}
\newcommand{\ot}{\otimes}
\begin{document}

\title[Dual Gabriel Theorem with applications]{Dual Gabriel Theorem with applications}
\author[Xiao-Wu Chen,  Hua-Lin Huang and Pu Zhang]{Xiao-Wu Chen,
Hua-Lin Huang  and Pu Zhang$^*$}
\thanks{$^*$ The corresponding author}
\thanks{Supported in part by the National Natural Science Foundation of China (Grant No. 10271113 and No. 10301033)
and the Europe Commission AsiaLink project  ``Algebras and
Representations in China and Europe$"$ ASI/B7-301/98/679-11}

%
\maketitle

\begin{center}
Department of Mathematics \\University of Science and Technology
of China \\Hefei 230026, Anhui, P. R. China
\end{center}

\vskip10pt

\begin{center}
Department of Mathematics\\ Shanghai Jiao Tong University\\
Shanghai 200030, P. R. China
\end{center}
\vskip10pt
\begin{center} xwchen$\symbol{64}$mail.ustc.edu.cn \ \ \ \
hualin$\symbol{64}$ustc.edu.cn

\vskip10pt

pzhang$\symbol{64}$ustc.edu.cn\end{center}

\vskip10pt


\begin{abstract}
 We introduce the quiver of a bicomodule over a cosemisimple coalgebra. Applying this to the
 coradical $C_0$ of an arbitrary coalgebra $C$, we give an alternative definition of the Gabriel quiver of
 $C$, and then show that it coincides with the known $\operatorname {Ext}$ quiver of
 $C$ and the link quiver of $C$.
 The dual Gabriel theorem for a coalgebra with separable
 coradical is  obtained, which generalizes the corresponding result for
a pointed coalgebra. We also give a new description of $C_1$ of any coalgebra $C$,
which can be regarded as a generalization of the first part of the well-known Taft-Wilson
Theorem for pointed coalgebras. As applications, we give a
characterization of locally finite coalgebras via their Gabriel quivers,
and a property of the Gabriel quiver of a quasi-coFrobenius coalgebra.

\end{abstract}

\vskip60pt

\section{Introduction and Preliminaries}

\vskip 10pt

\subsection{} In  the representation theory of finite-dimensional algebras,
 quiver is a fundamental tool.

A finite-dimensional algebra $A$ over a field $K$ is
called elementary, if the quotient algebra of $A$ modulo the Jacobson radical is
isomorphic to a product of $K$ as $K$-algebras, and basic if this quotient is  isomorphic to a product
of division $K$-algebras.  A theorem due
to Gabriel says that an elementary  $K$-algebra $A$ is isomorphic to the path algebra
$KQ(A)^a$  modulo an admissible ideal, where $Q(A)$ is the Gabriel quiver of $A$
(see Auslander-Reiten-Smal$\phi$ [ARS], Theorem 1.9,
 and Ringel [Rin], p.43).
Since any finite-dimensional algebra is Morita equivalent to a uniquely
determined basic algebra, and a basic algebra over an
algebraically closed field is elementary, it follows that any
finite-dimensional algebra $A$ over an algebraically closed field is
Morita equivalent to $KQ(A)^a$ modulo an admissible ideal. On the other
 hand, the Auslander-Reiten quiver of a finite-dimensional
 algebra $A$, which is defined by the indecomposable $A$-modules
 and  irreducible maps, is an essential approach and technique
 in studying the representations of $A$ (see e.g. [ARS] and [Rin]).

\vskip10pt

\subsection{} As pointed out by Chin and Montgomery in [CMon],
by the fundamental theorem for coalgebras
(i.e., every comodule is a sum of its finite-dimensional
subcomodules; in particular, every simple coalgebra is finite-dimensional),
it is reasonable to expect that the quiver technique for algebras could
be extended to coalgebras.

In fact, in the past  few years, there are several works towards
this direction. The path algebra construction has been
dualized by Chin and Montgomery to get a path coalgebra;
the $\operatorname {Ext}$ quiver of
coalgebra $C$ has been introduced and then a dual version of the Gabriel theorem
for coalgebras has been given in [CMon] (here $C$ is not necessarily finite-dimensional). Montgomery also
introduces the link quiver of coalgebra $C$ by using the wedge of
simple subcoalgebras of $C$. This link quiver is isomorphic to the
$\operatorname {Ext}$ quiver, up to multiple arrows, and it is connected if and only if
$C$ is an indecomposable coalgebra; using this she  proved
that a pointed Hopf algebra is a crossed product of a group algebra over
the indecomposable component of the identity element, see [M2]. On the other
hand, the almost split sequences and the Auslander-Reiten quivers
for coalgebras turn out to be also very useful in studying the
comodules of coalgebras. See Chin, Kleiner, and Quinn [CKQ],
and Simson [Sim].

\vskip10pt

There are also several works to construct neither commutative nor
cocommutative Hopf algebras via quivers. In [C] Cibils determined
all the graded Hopf structures with length grading on path
algebra $KZ_n^a$ of basic cycle $Z_n$. In [CR1]
Cibils and Rosso studied graded Hopf structures on path algebras.
In [GS] E. Green and Solberg studied Hopf structures on some
special quadratic quotients of path algebras. More recently,
Cibils and Rosso [CR2] introduced the notion of the Hopf quivers and
then classified all the  graded
Hopf algebras with length grading on path coalgebras.
Using the quiver technique all the monomial Hopf algebras
have been classified in [CHYZ], and in [WZ] a class of bi-Frobenius algebras
which are not Hopf algebras have been constructed via quivers.

\vskip 10pt

\subsection{} These quoted works inspire us to pay more attention to the quiver method towards
 coalgebras. Note that in algebra case the Gabriel quiver has an alternative definition rather than
the extensions of simple modules. In this paper, we first introduce the quiver of a
bicomodule over a cosemisimple coalgebra. Applying this to the
 $C_0$-$C_0$-bicomodule $C_1/C_0$, we give an alternative definition of the Gabriel quiver $Q(C)$ of
 $C$, where $C$ is an arbitrary coalgebra;
 and then show that $Q(C)$ coincides with the $\operatorname {Ext}$ quiver of
 $C$ introduced by Chin and Montgomery [CMon]. This is done in Section 2.

\vskip10pt

By definition a coalgebra $C$ is called pointed if each simple subcoalgebra of $C$ is of dimension one
(in finite-dimensional case, this  is exactly the dual of an elementary algebra), and basic if
the dual of each simple subcoalgebra of $C$ is a finite-dimensional division $K$-algebra.
As a dual of the result due to Gabriel as quoted in 1.1, Chin and Montgomery
proved that any pointed coalgebra is isomorphic to
a large subcoalgebra of the path coalgebra
of the $\operatorname {Ext}$ quiver of $C$ (for the notion of ``large subcoalgebra" see Remark 3.5).
Since any coalgebra is Morita-Takeuchi equivalent to a uniquely
determined basic coalgebra, and a basic coalgebra over an
algebraically closed field is pointed, it follows that any
coalgebra $C$ over an algebraically closed field is
Morita-Takeuchi equivalent to a large subcoalgebra of the path coalgebra
of the $\operatorname {Ext}$ quiver of $C$. See [CMon], Theorem 4.3.
In Section 3 (Theorem 3.1), we prove that
a coalgebra $C$ (over an arbitrary field $K$) with separable coradical $C_0$ is isomorphic to
a large subcoalgebra of the cotensor coalgebra ${\rm
Cot}_{C_0}(C_1/C_0)$. Note that  $C_0$ is always separable over an
algebraically closed field, and if $C$ is pointed then ${\rm Cot}_{C_0}(C_1/C_0)$ is
isomorphic to the path coalgebra $kQ^c$ of the Gabriel quiver $Q$
of $C$. In this way the dual of the Gabriel theorem for pointed coalgebras is extended to the one for coalgebras with
separable coradicals.

\vskip 10pt

For an arbitrary coalgebra $C$ with coradical $C_0 =
\oplus_{i \in I} D^i$, where $D^i$'s are simple subcoalgebras of
$C$, we show in Section 4 that there hold
$$C_1 = \sum_{i,j \in I} (D^i \wedge_C D^j)$$ and $$C_1/C_0\cong \oplus_{i,j \in I}(D^i \wedge_C D^j)/ (D^i + D^j).$$
See Theorem 4.1. This can be regarded as a generalization of the first part of the well-known Taft-Wilson
Theorem for pointed coalgebras, see Remark 4.2. As an application we unify the link quiver of a coalgebra
with the Gabriel quiver and the $\operatorname {Ext}$ quiver (Corollary 4.4).

\vskip 10pt

In the last two sections, we include two
applications of Theorem 3.1 and Theorem 4.1, by claiming that a coalgebra with separable coradical is
 locally finite if and only if its Gabriel quiver is locally
finite (Theorem 5.2); and the Gabriel quiver of
a non-simple quasi-coFrobenius coalgebra has no sources and no sinks (Theorem
6.2). In finite-dimensional case, Theorem 6.2 is dual
to the corresponding one for algebras.

\vskip10pt

\vskip10pt

In the following, all coalgebras and all tensor products are over a
fixed field $K$. For a
$K$-space $V$, denote the dual ${\rm Hom}_K(V,K)$ by $V^*$.

\vskip10pt

\subsection{}  Let $(C, \Delta_C, \varepsilon_C)$ be a coalgebra over
the ground field $K$,  where $\Delta_C$ and $\varepsilon_C$ are
the structure maps. A right $C$-comodule $(M,\rho)$ is a vector
space $M$ endowed with a structure map $\rho: M \longrightarrow M
\otimes C $ such that $(\rho \otimes Id) \circ \rho=(I \otimes
\Delta_C)\circ \rho$ and $(Id \otimes \varepsilon_C) \circ \rho = Id$,
where $Id$ denotes the identity map. Similarly one has left
$C$-comodules. Let $D$ be a coalgebra. By a $D$-$C$-bicomodule
$(M,\rho_l,\rho_r)$ we mean that $(M,\rho_l)$ is a left
$D$-comodule  and $(M,\rho_r)$ is a right $C$-comodule, satisfying
$(Id \otimes \rho_r)\circ \rho_l=(\rho_l \otimes Id)\circ \rho_r$.

\vskip 10pt

Let $(M, \rho)$ and $(N, \delta)$ be a right and a left
$C$-comodule, respectively. Then the cotensor product of $M$ and
$N$ over $C$ is defined to be the subspace of $M\otimes N$ given
by
\begin{center}
$M\square_C N= {\rm Ker} (\rho \otimes Id-Id \otimes \delta : M
\otimes N \longrightarrow M\otimes C \otimes N)$
\end{center}
If $M$ is a $D$-$C$-bicomodule and $N$ is a $C$-$D'$-bicomodule,
then $M \square_C N$ is a $D$-$D'$-bicmodule. The cotensor product
is associative, i.e., if in addition $L$ is a
${D'}$-${C'}$-bicomodule, then
 $(M \square_C N) \square_{D'} L \simeq M \square_C(N
 \square_{D'}L)$ as $D$-$C'$-bicomodules.

 \vskip 10pt

Let $(C,\Delta_C, \varepsilon_C)$ be a coalgebra, and $(M,\rho_l,
\rho_r) $  a $C$-$C$-bicomodule. Write $\rho_l(m)= \sum m_{-1}
\otimes m_0$ and $\rho_r(m)= \sum m_0 \otimes m_1 $ for every $m
\in M$. Define $M^{\square 0}=C$, $M^{\square 1}=M$ and
$M^{\square n}=(M^{\square n-1}) \square_C M$ for any $n \geq 2$.
Note that $M^{\square n}$ is a subspace of $M ^{\otimes n}$ for
all $n \geq 1$. If $\sum m^1\otimes \cdots \otimes m^n\in
M^{\square n}$, we write it as $\sum m^1 \square \cdots \square
m^n$. Define the cotensor coalgebra $ {\rm Cot}_C (M)$. As a
vector space, $\rm{Cot}_C(M)= \oplus_{ i=0 }^{\infty} M^{\square
i}$. The counit $\varepsilon$ is given by $\varepsilon
|_{M^{\square i}}=0 $ for $i \geq 1$ and $\varepsilon|_{M^{\square
0}}=\varepsilon_C$; the comultiplication $\Delta$ of $ {\rm Cot}_C
(M)$ is defined as $\Delta|_{M^{\square 0}}=\Delta_C$,
$\Delta(m)=\rho_l(m)+\rho_r(m)=\sum m_{-1} \otimes m_0 + m_0
\otimes m_1$ for all $m \in M$, and in general, if $\sum m^1
\square \cdots \square m^n \in M^{\square n}$ ($n \geq 1$), then
\begin{align*}
\Delta(\sum m^1 \square \cdots \square m^n)&=\sum (m^1)_{-1}
\otimes
( (m^1)_0 \square \cdots \square m^n)\\
&+ \sum_{i=1}^{n-1} (m^1 \square \cdots \square m^i) \otimes
(m^{i+1} \square \cdots \square m^n)\\
&+ \sum( m^1 \square \cdots \square (m^n)_0) \otimes (m^n)_1 \\ &
\in (C \otimes M^{\square n})\\ & \bigoplus_{i=1}^{n-1}(M^{\square
i} \otimes M^{\square
  (n-i)})\\
& \oplus M^{\square n }\otimes
 C
 \\
& \subseteq {\rm Cot}_C{M} \otimes {\rm Cot}_C(M).
\end{align*}
One can verify that $\Delta$ is well-defined and $({\rm
Cot}_C(M), \Delta, \varepsilon)$ is a coalgebra.

\vskip 10pt

\begin{rem}
In case $C$ is cosemisimple,  then coalgebra ${\rm Cot}_C(M)$
is coradically graded, i.e., $\{\oplus_{i \leq n} M^{\square i} |
n=0,1, \cdots\}$ is its coradical filtration (see [CMus, Sect.2]).
\end{rem}

\subsection{} By a quiver, we mean an oriented graph $Q =
(Q_0, Q_1, s,t)$ with $Q_0$ the set of vertices and $Q_1$ the set
of arrows, where $s,t$ are two maps from $Q_1$ to $Q_0$. For
$\alpha\in Q_1$, $s(\alpha)$ and $t(\alpha)$ denote the starting
and terminating vertex of $\alpha$, respectively. Note that
quivers considered here could be infinite. \vskip 10pt \par

Recall that the  path coalgebra  $KQ^c$ of a qiver $Q$ is defined
as  follows (see [CMon]). As a vector space $KQ^c$ has a basis
consisting of paths in $Q$; the comultiplication  is given by
\begin{align*}
\Delta(p) =  \alpha_l \cdots \alpha_1\otimes s(\alpha_1)+
\sum_{i=1}^{l-1}\alpha_l \cdots \alpha_{i+1} \otimes
\alpha_{i}\cdots \alpha_1+t(\alpha_l)\otimes\alpha_l \cdots
\alpha_1
\end{align*}
for each path $p = \alpha_l\cdots \alpha_1$ with each $\alpha_i\in
Q_1$; and $\varepsilon(p)=0$ if $l\ge 1$, and $1$ if
$l =0$. Then $KQ^c$ is a pointed coalgebra with coradical
filtration $C_n=kQ_0 \oplus \cdots \oplus kQ_n$, where $KQ_n$ is
the $K$-space with basis the set of all paths of length $n$.
\vskip 10pt
\par

\begin{rem}
Note that a path coalgebra is a special case of a cotensor
coalgebra: for every quiver $Q$, $KQ^c \simeq {\rm
Cot}_{kQ_0}(KQ_1)$ as coalgebras, where the bicomodule structure
of $KQ_1$ is given as
 $\rho_l(\alpha):= t(\alpha)
\otimes \alpha $  and $ \rho_r(\alpha):= \alpha \otimes s(\alpha)$
for each $\alpha \in Q_1$.

On the other hand, if $C$ is a pointed coalgebra, then $C_1/C_0$
becomes a $C_0$-$C_0$-bicomodule and $KQ(C)^c \simeq {\rm
Cot}_{C_0}(C_1/C_0)$, where $Q(C)$ is the Gabriel quiver of $C$,
as defined in 2.2 below, see also Remark 3.5 below.
\end{rem}

\vskip10pt

\section{The Gabriel quiver of a coalgebra}

\vskip 10pt

\subsection{} Quiver of a bicomodule over a cosemisimple coalgebra

\vskip10pt

Let $D$ be a cosemisimple $K$-coalgebra and let $(M, \rho_l,
\rho_r)$ be a $D$-$D$-bicomodule.  We associate a quiver with the
bicomodule $M$.

Write $D =\oplus_{i \in I} D^i$, where each $D^i$ is a simple
subcoalgebra of $D$. Set
$$^iM^j =\{\ m \in M \ |\
\rho_l(m) \in D^i \otimes M, \ \ \ \ \rho_r(m) \in M \otimes
D^j\}$$ for each $i,j \in I$. Then $M=\oplus_{i,j \in I}$ $^i M
^j$, and each $^iM^j$ is naturally a $D^i$-$D^j$-bicomodule, and
hence a $(D^j)^*$-$(D^i)^*$-bimodule. Since each $(D^i)^*$ is a
simple algebra, it follows that $(D^i)^* \simeq
M_{n_i}(\Delta_i)$, where $\Delta_i$ is a finite-dimensional
division algebra over $K$. For each $(D^i)^*$ fix a primitive
idempotent $e_i$. Set
\begin{center}
$t_{ij} ={\rm dim}_{K} (e_i. ^jM^i . e_j)$
\end{center}
for each pair of $i,j \in I$, where the dots denote the module
action. Note that $t_{ij}$ is independent of the choice of
$e_i$'s, because they are mutually conjugate in
$M_{n_i}(\Delta_i)$.

Define the quiver $Q(D,M)$ of a $D$-$D$-bicomodule $M$ as follows:
the set of vertices is $I$, and for any $i, j \in I$, the number
of arrows from $i$ to $j$ is $t_{ij}$.

\vskip 10pt

\begin{rem} 1. We admit the case that $I$ is an infinite set and $t_{ij}$
is  infinite, i.e., the quiver $Q(D,M)$ could be infinite.

\vskip10pt

2. If $D$ is group-like (i.e. the set of group-like elements $G =
G(D)$ forms a basis for $D$, or equivalently, $D$ is cosemisimple
pointed), then the quiver $Q(D,M)$ is simply interpreted as
follows. The set of vertices is $G$, and for any $g, h \in G$, the
number of arrows from $g$ to $h$ is $t_{gh}$, where $t_{gh}= \dim_K
{^h M^g}$, and $^hM^g=\{m \in M \ |\ \rho_l(m)= h \otimes m, \ \
\rho_r(m)=m \otimes g\}$.
\end{rem}

 \vskip10pt

\subsection{} The Gabriel quiver of a coalgebra

\vskip10pt

Let $(C,\Delta)$ be a coalgebra with coradical filtration
$\{C_n\}$.  Set $\pi_0: C \longrightarrow C/{C_0}$ to be the
canonical projection. Define map $\tilde{\rho_l}: \
C\longrightarrow C\otimes C/C_0$ by
 $\tilde{\rho_l}=(Id
\otimes \pi_0)\circ \Delta$, and map $\tilde{\rho_r}: \
C\longrightarrow C/C_0\otimes C$ by $\tilde{\rho_r}=(\pi_0 \otimes
Id)\circ \Delta$. Since $\tilde{\rho_l}(C_0)=0$,
$\tilde{\rho_r}(C_0)=0$, and $\Delta(C_1) \subseteq C_0 \otimes
C_1 + C_1 \otimes C_0 $, it follows that $\tilde{\rho_l}$ and
$\tilde{\rho_r}$ induce two maps $\rho_l: C_1/C_0 \longrightarrow
C_0 \otimes C_1/C_0$, and $\rho_r: C_1/C_0 \longrightarrow C_1/C_0
\otimes C_0$, respectively. It is clear that $(C_1/C_0, \rho_l,
\rho_r)$ is a $C_0$-$C_0$-bicomodule.

\vskip10pt

\begin{defn} The Gabriel quiver $Q(C)$ of a
coalgebra $C$ is defined to be the quiver $Q( C_0, C_1/C_0)$ of
$C_0$-$C_0$-bicomodule $C_1/C_0$.

More precisely, let $C_0 = \oplus_{i \in I} D^i$, where $D^i$'s are
simple subcoalgebras, and $e_i$ be a fixed  primitive idempotent
of $(D^i)^*$. Then vertices of $Q( C_0, C_1/C_0)$ are $i\in I$,
and there are $t_{ij}= {\rm dim}_{K} e_i . ^j( C_1/C_0)^i . e_j$
arrows from $i$ to $j$.
\end{defn}

\vskip10pt

\subsection{} The $\operatorname {Ext}$ quiver of a coalgebra

\vskip10pt

Let $C$ be a $K$-coalgebra. Recall the definition of the
$\operatorname {Ext}$ quiver of $C$ introduced by Chin and
Montgomery ([CMon], p.45). Let $\{S_i \ | \ i \in I\}$ be the
complete set of isoclasses of right simple $C$-comodules.  The Ext
quiver of $C$ is an oriented graph with vertices indexed by $I$, and
there are ${\rm dim}_{K}{\rm Ext}^1 (S_i, S_j)$ arrows from $i$ to
$j$ for any $i,\ j \in I$.

\vskip10pt

Note that in Simson [Sim], p.468, the $\operatorname {Ext}$ quiver
is also called the Gabriel quiver. In [M2] Montgomery also
introduced the link quiver of coalgebra $C$ by using the wedge of
simple subcoalgebras of $C$. This link quiver is isomorphic to the
$\operatorname {Ext}$ quiver, up to multiple arrows, see [M2,
Theorem 1.7] . Moreover, it is shown in [M2, Corollary 2.2 ]  that
$\operatorname {Ext}$ quiver of $C$ is connected if and only if
$C$ is an indecomposable coalgebra.

\vskip10pt The main result of this section is

\vskip10pt

\begin{thm} The Gabriel quiver of $C$ coincides with the Ext quiver of $C$.
\end{thm}

\vskip 10pt

\subsection{} In order to prove the result, we need some preparations.

For a right $C$-comodule $M$, denote by $E(M)$ its injective hull,
which always exists and $soc(M)
= soc(E(M))$ (see [G] or [DNR], Chap.2). Let $\{S_i\ |\ i \in I\}$ be the complete set of
isoclasses of simple right $C$-comodules.  Then as right
$C$-comodule we have $D^i \simeq n_i S_i$ and $C \simeq \oplus_{i
\in I} E(D^i)\simeq \oplus_{i \in I} n_i E(S_i)$. Note that
$(D^i)^* \simeq M_{n_i} (\Delta_i)$, where $\Delta_i$ is a
finite-dimensional division algebra over $K$.

\vskip 10pt

\begin{lem} Suppose ${\rm dim}_K \Delta_i= d_i$ for each $i \in I$. We
have
$$ soc (E(D^i)/D^i)= \oplus_{j \in I} \frac { n_i t_{ji}}{d_j}
S_j$$ and $$soc (E(S_i)/S_i) \simeq  \oplus_{j \in I}
\frac{t_{ji}}{d_j} S_j.$$
\end{lem}
\vskip 5pt \noindent{\bf Proof} \quad  Recall that  $(C_1/C_0,
\rho_l, \rho_r)$ is  a $C_0$-$C_0$-bicomodule. Set

\begin{center}
$ ^i(C_1/C_0)=\{x \in C_1/C_0 \ |\ \rho_l(x) \in D^i\otimes
(C_1/C_0) \}$.
\end{center}
Then $ ^i(C_1/C_0)$ is a left $D^i$-comodule.
\par \vskip 5pt

We may identify $C$ with $\oplus_{i \in I} E(D^i)$. It follows
that as a right $C$-comodule we have
$$C_1 = \oplus_{i \in I}(E(D^i)\cap C_1).$$

(In fact, for each $c \in C_1$, we have $c= \sum c_i$ with $c_i
\in E (D^i)$. Since $\Delta(c_i) \in E(D^i) \otimes C$, it follows
that $\Delta(c_i)=\sum_j d_{ij}\otimes c_{ij}$ with $d_{ij}\in
E(D^i)$ and $\{d_{ij}\}$ linearly independent for each $i\in I$.
Then $\Delta(c)= \sum_{i,j} d_{ij}\otimes c_{ij} \in C_1 \otimes
C_1$. Since $\{d_{ij}\}$ is linearly independent, it follows that
each $c_{ij}$ is contained in $C_1$, and hence by the counitary
property, each $c_i \in C_1$.)

Thus
$$C_1/C_0 = \oplus_{i \in I}(E(D^i)\cap C_1)/D^i$$
as right $C_0$-comodules, and hence $(E(D^i)\cap C_1)/D^i$ is a
cosemisimple $C_0$-comodule, which implies $(E(D^i)\cap
C_1)/D^i\subseteq soc(E(D^i)/D^i).$ While $soc(C/C_0) =C_1/C_0$
(see [M1], p.64), it follows that

\begin{align*}soc(C/C_0) &= soc((\oplus_{i \in I}E(D^i))/(\oplus_{i \in I}D^i)) = \oplus_{i \in I}
soc(E(D^i)/D^i)\\ & = C_1/C_0 = \oplus_{i \in I}(E(D^i)\cap
C_1)/D^i.\end{align*}
This forces $$soc(E(D^i)/D^i) = (E(D^i)\cap
C_1)/D^i.$$

We claim $soc(E(D^i)/D^i)\subseteq {^i(C_1/C_0)}$, and then by $C_1/C_0
=\oplus_{i \in I} { ^i(C_1/C_0)}$ we have
$$soc(E(D^i)/D^i)= {^i(C_1/C_0)}.$$

To see this, note that
$$\tilde{\rho_l}((E(D^i))\subseteq E(D^i)\otimes C, \ \
\tilde{\rho_l}(C_1)\subseteq C_0\otimes C_1,$$ it follows that
$$\tilde{\rho_l}((E(D^i)\cap C_1)\subseteq (E(D^i)\cap C_0)\otimes
C_1.$$ Note that $E(D^i)\cap C_0\subseteq soc(E(D^i)) = D^i$, it
follows that

$$\tilde{\rho_l}(E(D^i)\cap C_1) \subseteq D^i\otimes C_1$$ and hence

$$\rho_l(soc(E(D^i)/D^i)) = \rho_l((E(D^i)\cap C_1)/D^i)\subseteq D^i\otimes (C_1/C_0). \ \ $$
That is $soc(E(D^i)/D^i)\subseteq {^i(C_1/C_0)}$.
This proves the assertion.
\vskip10pt

Note that $^i (C_1/C_0)^j$ is a $D^i-D^j$-bicomodule, and hence a
$(D^j)^*-(D^i)^*$-bimodule. Thus $^i (C_1/C_0)^j.e_i$ is a left
$(D^j)^*$-module, and hence a right $D^j$-comodule, where $e_i$ is
a primitive element of $(D^i)^*$. Thus we have
$${^i (C_1/C_0)^j}.e_i=
m_jS_j$$ as a right $D^j$-comodule, for some non-negative integer
$m_j$. Since $t_{ji} ={\rm dim}_K (e_j. ^i(C_1/{C_0})^j .e_i)$ and ${\rm dim}_K S_j = n_j d_j$, it follows that
$$m_jn_jd_j  = {\rm dim}_K
{^i (C_1/C_0)^j.e_i} = {\rm dim}_K n_j (e_j. ^i(C_1/{C_0})^j .e_i)
= n_j t_{ji}$$ and hence $m_j = \frac{t_{ji}}{d_j}$. It follows
that

\begin{align*}
soc(E(D^i)/D^i) &= {^i (C_1/C_0) } \\ &= \oplus_{j \in I} {^i (C_1/C_0)^j} \\
                     &= \oplus_{j \in I} n_i \ {^i (C_1/C_0)^je_i}\\
                     & = \oplus_{j \in I}\frac{n_i t_{ji}} {d_j} S_j.\ \ \ \ \ \ \square
\end{align*}

\vskip 10pt

\subsection{}Proof of Theorem 2.3

Since $D^j \simeq n_j S_j$, it suffices to compute ${\rm
Ext}^{1}(S_i, D^j)$. For this, we take a minimal injective
resolution of $D^j$ (see [D])
$$ 0 \longrightarrow D^j \stackrel {d_0}{\longrightarrow} E_0 \stackrel {d_1}{\longrightarrow}
E_1 \stackrel{d_2}{\longrightarrow} E_2 \longrightarrow \cdots
$$
 where $E_0=E(D^j)$ and $E_1=E(E_0/D^j)$. Since
$$Im(d_0) = soc(E_0), \ \ Im(d_1) \supseteq soc(E_1),$$
it follows that for every comodule  map $g: S_i \longrightarrow
E_0$ we have
 $ d_1 \circ g=0$, and that for every
 comodule map $f: S_i \longrightarrow E_1$ we have $d_2 \circ
 f=0$.   It follows that
 $${\rm Ext}^1(S_i,D^j)={\rm Hom}_{C} (S_i, E_1)= {\rm Hom}_{C}(S_i, soc (E_0/D^j)),$$
 here we have used the  fact  $soc(E_1)=soc (E_0/D^j)$. By Lemma 2.3 we have  $soc(E_0/D^j)= \oplus_{i \in I}
 \frac{n_jt_{ij}}{d_i}
 S_i$.   Combining this with the fact ${\rm Hom}_{C} (S_i,
 S_i)=\Delta_i$,
  we have $$t_{ij}=  \frac{1}{n_j}{\rm dim}_{K} Ext^1
 (S_i,D^j)= {\rm dim}_K Ext^1(S_i, S_j).$$ This completes the proof.  \hfill $\square$

 \vskip 10pt

\begin{rem}  Recall that two coalgebras $C$ and $D$ are said to be
Morita-Takeuchi equivalent, if the categories of $C$-comodules and
$D$-comodules are equivalent (see [T]). Then by Theorem 2.3 two
coalgebras have the same Gabriel quiver provided they are
Morita-Takeuchi equivalent. \end{rem} \vskip10pt

\section{The dual Gabriel Theorem}

\vskip10pt

\subsection{} Let $L$ be a field extension of $K$, and $C$ a $K$-coalgebra.
Then $C \otimes L$ is naturally an $L$-coalgebra. A coalgebra is
called separable provided that $C \otimes L$ is cosemisimple for
any field extension $L$. Note that $C$ is separable if and only if
$C\otimes C^{cop}$ is cosemisimple.  (In fact, cosemisimple
coalgebras are direct sum of simple coalgebras, thus this follows
by dualizing [DK], Theorem 6.1.2.) For example, a group-like
coalgebra $C$ is separable. If $K$ is of characteristic zero, then
any cosemisimple coalgebra is separable. Note that the coradical
$C_0$ of $C$ is separable if $C$ is pointed or if $K$ is
algebraically closed.

\vskip10pt

The main result of this section is

\begin{thm}
Let $C$ be a coalgebra with separable coradical $C_0$. Then there
exists a coalgebra embedding $i: \ C \hookrightarrow {\rm
Cot}_{C_0}(C_1/C_0)$ with $i(C_1) = C_0\oplus C_1/C_0$.
\end{thm}

\vskip 10pt

\subsection{} To prove Theorem 3.1, one needs the
following fundamental lemma, which gives the universal mapping
property of cotensor coalgebras.

\vskip 10pt

Let $C$ and $D$ be coalgebras and $f: D \longrightarrow C$ be a
coalgebra map. Then $D$ becomes a $C$-$C$-bicomodule via $f$: the
left and right comodule structure maps are given by $(f \otimes
Id)\circ \Delta_D$ and $(Id \otimes f) \circ \Delta_D$,
respectively.

\par \vskip 5pt

\begin{lem}
Let $C$ and $D$ be coalgebras and $M$  a $C$-$C$-bicomodule. Given
a coalgebra map $f_0: D \longrightarrow C$,  and a $C$-$C$-
bicomodule map $f_1: D \longrightarrow M$ with property that $f_1$
vanishes on the coradical $D_0$ of $D$, where the $C$-$C$-
bicomodule structure $D$ is given via $f_0$. Then there exists a
unique coalgebra map
$$F: D \longrightarrow \rm{Cot}_C(M)$$ with $\pi_i \circ F=f_i$ $(i=0,1)$, where each $\pi_i
:\rm{Cot}_C(M) \longrightarrow M^{\square i}$ is the canonical
projection.
\end{lem}

\noindent {\bf Proof} \quad  Set $\Delta_0$ to be the identity map
of $D$, $\Delta_1= \Delta_D$ and  define $\Delta_{n+1}=(\Delta_D
\otimes Id) \circ \Delta_n$ for all $n \geq 1$, where $Id$ denotes
the identity map of $D^{\otimes n}$.  It is easy check that
$\Delta_n(D) \subseteq D^{\square n+1}$ and    $f_1 ^{\otimes n+1}
\circ \Delta_n(D) \subseteq M^{\square n+1}$ for each $n \geq 1$.
We claim that  $F: D \longrightarrow {\rm Cot}_C(M)$ given by
\begin{center}
$F(d)=f_0(d)+ \sum_{n=0}^{\infty} f_1^{\otimes {n+1}}\circ
\Delta_n(d)$
\end{center}
for each $d \in D$, is well-defined.
\par \vskip 5pt

In fact, we have $\bigcup_{n\ge 0} D_n=D$, where $\{D_n\}$ is the
coradical filtration of $D$ (see e.g. [M1], Theorem 5.2.2). Thus
for each $d \in D_n$, $f_1^{\otimes m+1} \Delta_m (d)=0$ for all
$m \geq n$. This is because
$$\Delta_m(D_n) \subseteq \sum_{i_0+i_1+\cdots + i_m = n} D_{i_0} \otimes \cdots \otimes D_{i_m}$$
and $f_1$ vanishes on $D_0$.  Thus $F$ is well-defined. Moreover,
$F$ is a coalgebra map with $\pi_i \circ F=f_i$ $(i=0,1)$.
\par \vskip 5pt

It remains to prove the uniqueness of coalgebra map  $F: D
\longrightarrow {\rm Cot}_C(M)$ with $\pi_i \circ F=f_i$
$(i=0,1)$. Set $f_n = \pi_n \circ F$ for each $n$. It suffices to
prove that $f_n= f_1^{\otimes n} \circ \Delta_{n-1} $ for every $n
\geq 1$. Use induction on $n$. Assume that $f_m= f_1^{\otimes m}
\circ \Delta_{m-1}$, $m\ge 1$. Consider $f_{m+1}$. For every $d
\in D$, write $\Delta_D(d)= \sum d_1 \otimes d_2$. Since $F$ is a
coalgebra map, it follows that $\Delta (F(d))= (F\otimes
F)\Delta_D(d)$. Writing out the both sides explicitly we have

$$\Delta (F(d))= \sum_n \Delta (f_n(d))$$
with $$\Delta (f_n(d))\in C \otimes M^{\square n}\oplus M \otimes
M^{\square (n-1)} \oplus \cdots \oplus M^{\square (n-1)} \otimes M
\oplus M^{\square n} \otimes C;$$ and

\begin{align*}(F\otimes F)\Delta_D(d) &= \sum_{(d)} F(d_1)\otimes F(d_2)\\ &=
\sum_{n}\sum_{(d), i+j=n} f_i(d_1) \otimes f_j(d_2)\end{align*}
with
\begin{align*}&\sum_{(d), i+j=n} f_i(d_1) \otimes
f_j(d_2)\\ &\in C \otimes M^{\square n}\oplus M \otimes M^{\square
(n-1)} \oplus \cdots \oplus M^{\square (n-1)} \otimes M \oplus
M^{\square n} \otimes C.\end{align*} It follows that $$\Delta (f_{n} (d))=
\sum_{(d), i+j=n} f_i(d_1) \otimes f_j(d_2), \ \ \forall n\ge 2.$$
Note that $f_n(d)\in M^{\square n}$, and $f_i(d_1) \otimes
f_j(d_2)\in M^{\square i}\otimes M^{\square j}$. By the definition
of the comultiplication $\Delta$ of ${\rm Cot}_C(M)$ and by
comparing the terms belonging to $M^{\square i}\otimes M^{\square
j}$ with $i\ne 0\ne j$ and $i+j = n$, we obtain
$$ f_{n}(d)=\sum_{(d)}f_{i} (d_1) \otimes f_j(d_2).$$
In particular we have by induction
\begin{align*} f_{m+1}(d)=&\sum_{(d)}f_{m} (d_1) \otimes f_1(d_2) \\
          =& \sum f_1^{\otimes m } \circ \Delta_{m-1} (d_1) \otimes f_1(d_2)
          \\
          =& f_1^{\otimes {m+1}} \circ \Delta_m (d).
\end{align*}
This completes the proof. \hfill $\square$

\vskip 10pt

 To complete the proof of Theorem 3.1, we also need the dual Wedderburn-Malcev
 theorem (see [A], Theorem 2.3.11, or [M1], Theorem 5.4.2)
 and another lemma due to Heyneman-Radford (see [HR], or [M1], Theorem 5.3.1).

 \vskip5pt

\begin{lem}(Dual Wedderburn-Malcev theorem)
 Let $C$ be a coalgebra
with separable coradical. Then there is a coideal $I$ such that
$C= C_0\oplus I $, i.e, there is a coalgebra projection $\pi: C
\longrightarrow C_0$ such that $\pi|_{C_0}=Id$.
\end{lem}
 \vskip 5pt
\begin{lem} (Heyneman-Radford)
 Let $C$ and $D$ be coalgebras and $f: C \longrightarrow D$   a
coalgebra map. Then $f$ is injective if and only if $f|_{C_1}$ is
injective.
\end{lem}
\vskip 5pt

\subsection{} Now we are in a position to prove Theorem 3.1.

By the dual Wedderburn-Malcev theorem, there is a coideal $I$ of
$C$ such that $C=C_0 \oplus I$. Thus we have  a coalgebra
projection $f_0: C \longrightarrow C_0$ such that $f_0|_{C_0}=Id$.
Note that $C$ becomes a $C_0$-$C_0$-bicomodule via $f_0$, $I$ is a
$C_0$-$C_0$-subbicomodule of $C$. Set $C_{(1)}=C_1\cap I$. Then
$C_1= C_0 \oplus C_{(1)}$. Note that $C_{(1)}$ is a
$C_0$-$C_0$-subbicomodule of $I$ and the canonical vector space
isomorphism $\theta: C_{(1)} \simeq C_1/C_0$ is a
$C_0$-$C_0$-bicomodule map.

\par \vskip 5pt

View $I$ as a right $C_0\otimes {C_0}^{cop}$-comodule and
$C_{(1)}$ its subcomodule. Since $C_0$ is separable, it follows
that there exists a $C_0 \otimes {C_0}^{cop}$-comodule
decomposition $I=C_{(1)} \oplus J$. Thus we have a
$C_0$-$C_0$-bicomodule projection $p: I \longrightarrow C_{(1)}$
such that $p|_{C_{(1)}} = Id$. Define a map $f_1= \theta \circ p
\circ f_0^{'}$ from $C$ to $C_1/C_0$ where $f_0^{'}: C
\longrightarrow I$ is
  the canonical projection. Clearly $f_1: C \longrightarrow
  C_1/C_0$ is a $C_0$-$C_0$-bicomodule map vanishing on $C_0$.
 Thus, by Lemma 3.3 we obtain a unique coalgebra map $i : C \longrightarrow {\rm
  Cot}_{C_0}(C_1/C_0)$ such that $\pi_0 \circ i =f_0$ and $\pi_1 \circ
  i=f_1$. Clearly $i(C_1) = C_0 \oplus  C_1/C_0$.
By Lemma 3.5, $i$ is  injective.  This completes the proof. \hfill
$\square$ \par

\vskip10pt

\begin{rem} 1. Note that if $C$ is pointed, then ${\rm Cot}_{C_0}(C_1/C_0)$ is
isomorphic to the path coalgebra $kQ^c$ of the Gabriel quiver $Q$
of $C$.

(In order to see this, just note that $KQ^c$ and ${\rm
Cot}_{C_0}(C_1/C_0)$ both have the universal mapping property, and
then the assertion follows from Lemma 3.4.)

It follows from the result above that a pointed coalgebra embeds
in the path coalgebra of the Gabriel quiver of $C$. This has been
obtained by Chin and Montgomery in [CMon], Theorem 4.3. See also
[Rad], Corollary 1.

\vskip10pt

2. Recall that a subcoalgebra $D$ of a coalgebra $C$ is said to be
large provided that $D$ contains $C_1$. By the definition of the
Gabriel quiver, a large subcoalgebra $D$ of $C$ has the same
Gabriel quiver as $C$. Then Theorem 3.1 says that a coalgebra $C$
(not necessarily finite-dimensioal) with separable coradical is
isomorphic to a large subcoalgebra of cotensor coalgebra ${\rm
Cot}_{C_0}(C_1/C_0)$.

Recall that any finite-dimensional elementary algebra $A$ is
isomorphic to the path algebra of the Gabriel quiver of $A$ modulo
an admissible ideal (see e.g. [ARS], Theorem 1.9, or [Rin], p.43).
Thus, Theorem 3.1 can be regarded as a generalization of the dual
of this basic result for algebras (note the condition ``large" in
Theorem 3.1 just corresponds to the condition ``admissible" in the
case for algebras).
\end{rem}

\vskip 10pt

\section{A description of $C_1$}

\vskip10pt

\subsection{} Let $C$ be a coalgebra. Following [S], the wedge
of two subspaces $V$ and $W$ of $C$ is defined to be the subspace

$$V \wedge_C W: = \{c\in C\ | \ \Delta_C(c)\in V \otimes C + C \otimes W\}.$$

Let $C_0$ be the coradical of $C_0$, i.e., $C_0$ is the sum of of all simple subcoalgebras of $C$.
Recall that by definition
$C_n = C_0\wedge_C C_{n-1}$ for $n\ge 1$, and $\{C_n\}$ is called the coradical filtration of $C$.
Then $C_n$ is a subcoalgebra of $C$ with $C_n\subseteq C_{n+1}$, $C = \cup_{n\ge 0}C_n$, and
$\D C_n \subseteq\sum_{0\le i\le n} C_i\ot C_{n-i}$ (see e.g. [M1], 5.2.2).
For properties of wedges see [S], Chap. 9, and [HR], Section2).

\vskip 10pt

\subsection{} Let $(C, \Delta, \varepsilon)$ be a coalgebra with  dual
algebra $C^*$. For  $c \in C$ and $f \in C^*$, define

$$f \rightharpoonup c = \sum c_1 f(c_2)$$ and $$c \leftharpoonup f=
\sum f(c_1) c_2$$ where $\Delta(c)= \sum c_1 \otimes
  c_2$. Then it is well-known that
$C$ becomes a $C^*$-$C^*$-bimodule with $c = \varepsilon
\rightharpoonup c = c \leftharpoonup \varepsilon$ (see e.g. [M1],
1.6.5).

\vskip10pt

The following result gives a new description of $C_1$. We will use it in the next section, but also it
seems to be of independent interest.

\vskip10pt

\begin{thm} Let $C$ be a coalgebra with coradical $C_0 =
\oplus_{i \in I} D^i$, where $D^i$ are simple subcoalgebras of
$C$. Then

\vskip5pt\vskip5pt

(i) \ \ $C_1 = \sum_{i,j \in I} (D^i \wedge_C D^j).$

\vskip5pt\vskip5pt

(ii) \ \ $(D^i \wedge_C
D^j) \cap C_0= D^i + D^j, \  \ \forall i, j\in I.$

\vskip10pt

(iii)\ \ $C_1/C_0\cong \oplus_{i,j \in I}(D^i \wedge_C D^j)/ (D^i + D^j).$

\vskip5pt

\vskip5pt

(iv) \ \ $(D^i \wedge_C D^j)/ (D^i + D^j)\cong\ {^i(C_1/C_0)^j}, \  \ \forall i, j\in I.$
\end{thm}

\vskip10pt
\begin{rem} Recall that the set of group-like elements of a coalgebra $C$ is $G(C): = \{
\ 0\ne c\in C\ | \ \Delta(c) = c\otimes c\ \}$, and that a coalgebra $C$ is
said to be pointed if each simple subcoalgebra of $C$ is of
dimension one. Note that $C$ is pointed if and only if $C_0 = KG(C)$.
For $g,h \in G(C)$, denote by $P_{g,h}(C): = \{c
\in C \ |\ \Delta(c)= c\otimes g + h \otimes c\}$, the set of
$g,h$-primitive elements in $C$.  A $g,h$-primitive element $c$
is said to be non-trivial if $c \notin K (g-h)$.

Let $P'_{g,h}(C)$ be a subspace of $P_{g,h}(C)$ such that $$P_{g,h}(C)=P'_{g,h} (C) \oplus K(g-h).$$
Then the first part of the Taft-Wilson Theorem for pointed coalgebras says that
if $C$ is pointed, then
$$C_1 = KG(C)\oplus (\bigoplus_{g, h}P'_{g,h}(C)),$$
and hence
$$C_1/C_0 = \bigoplus_{g, h}P'_{g,h}(C) = \bigoplus_{g, h}(Kh \wedge_C Kg) / (Kg + Kh).$$
(For the last equality see e.g. [OZ], Lemma 4.2.)
From this point of view
Theorem 4.1 (iii) can be regarded as a form of the first part of the Taft-Wilson Theorem in general case.
\end{rem}

\vskip10pt

\subsection{}\noindent{\bf Proof of Theorem 4.1}

\vskip10pt

(i) \ \ On one hand, we have
$$C_1 = C_0\wedge_C  C_0 = (\sum_{i\in I} D^i)\wedge_C(\sum_{j \in I} D^j)
\supseteq \sum_{i,j \in I} (D^i \wedge_C D^j).$$

On the other hand, by an elementary argument in linear algebra, we
can write $C=C_0 \oplus V$ with a subspace $V$ such that
$\varepsilon(V)=0$. Take $\varepsilon_i \in C^*$ such that
$$\varepsilon_i|_{D^i}=\varepsilon, \ \ \ \varepsilon_i|_{D^j \oplus
V}=0 \ \ \ (j\neq i).$$ Then
$$\varepsilon(c) = \sum_{i\in I} \varepsilon_i(c), \ \ \ \forall \ c\in C,$$
and hence by the counitary property we have
$$c=  \sum_{i, j\in I}
(\varepsilon_j \rightharpoonup c \leftharpoonup \varepsilon_i),\ \
\ \forall \ c\in C.$$ While for $c \in C_1$ we claim
$$\varepsilon_j \rightharpoonup c \leftharpoonup \varepsilon_i \in
D^i \wedge_C D^j,$$ and then the assertion follows.

In order to see the claim, for $c\in C_1$ consider $\D^3(c)= (\D\ot Id \ot Id)(Id\ot\D)\D(c)$. For simplicity
we omit the sum in the following
\begin{align*}\D(c)& =  c_1\ot c_2\in
\sum_s C_1\ot D^s + \sum_t D^t\ot C_1;\\ \\
\D^2(c)&= (Id\ot\D)\D(c)
= c_{1}\ot c_{21}\ot c_{22}\\ & \in
\sum_s C_1\ot D^s\ot D^s + \sum_{t,k} D^t\ot C_1\ot D^k + \sum_{t, k} D^t\ot D^k\ot C_1;\\
\\
\D^3(c)&= (\D\ot Id \ot Id)(Id\ot\D)\D(c)
= c_{11}\ot c_{12}\ot c_{21}\ot c_{22}\\ & \in
\sum_s C_1\ot C_0\ot D^s\ot D^s
+
\sum_s C_0\ot C_1\ot D^s\ot D^s
 \\ & +
 \sum_{t, k} D^t\ot D^t\ot C_1\ot D^k
 +
 \sum_{t, k} D^t\ot D^t\ot D^k\ot C_1.
\end{align*}
By definition we have
$$\D(\varepsilon_j \rightharpoonup c \leftharpoonup \varepsilon_i)
= \varepsilon_i(c_{11})\varepsilon_j(c_{22})c_{12}\ot c_{21}.$$
If $$c_{11}\ot c_{12}\ot c_{21}\ot c_{22}\in \sum_s C_1\ot C_0\ot D^s\ot D^s,$$
then
$$\varepsilon_i(c_{11})\varepsilon_j(c_{22})c_{12}\ot c_{21}\in C_0\ot D^j;$$
if $$c_{11}\ot c_{12}\ot c_{21}\ot c_{22}\in \sum_s C_0\ot C_1\ot D^s\ot D^s,$$
then
$$\varepsilon_i(c_{11})\varepsilon_j(c_{22})c_{12}\ot c_{21}\in C_1\ot D^j;$$
if $c_{11}\ot c_{12}\ot c_{21}\ot c_{22}\in \sum_{t, k} D^t\ot D^t\ot C_1\ot D^k$,
then
$$\varepsilon_i(c_{11})\varepsilon_j(c_{22})c_{12}\ot c_{21}\in D^i\ot C_1;$$
if $c_{11}\ot c_{12}\ot c_{21}\ot c_{22}\in \sum_{t, k} D^t\ot D^t\ot D^k\ot C_1,$
then
$$\varepsilon_i(c_{11})\varepsilon_j(c_{22})c_{12}\ot c_{21}\in D^i\ot C_0.$$
Thus, in all the cases we have
$$\varepsilon_i(c_{11})\varepsilon_j(c_{22})c_{12}\ot c_{21}\in D^i\ot C + C\ot D^j.$$
This proves $\varepsilon_j \rightharpoonup c \leftharpoonup \varepsilon_i\in D^i\wedge_C D^j$.

\vskip10pt

(ii) \ \ This is straightforward (or, follows from [HR], Lemma
2.3.1).

\vskip10pt

(iii) and (iv) \ \ Since $(D^i \wedge_C D^j) \cap C_0= D^i + D^j$ by (ii), it follows that there is
a coalgebra embedding
\begin{center} $(D^i \wedge_C D^j)/ (D^i + D^j)\cong ((D^i \wedge_C D^j)+C_0)/ C_0
\hookrightarrow C_1/C_0$.
\end{center}
 By the construction of
$C_0$-$C_0$-bicomodule structure maps $\rho_l$ and $\rho_r$ of
$C_1/C_0$, one observes that  $$(D^i \wedge_C D^j)/ (D^i + D^j)
\hookrightarrow {^i (C_1/C_0) ^j}.$$
It follows from (i) that
\begin{align*}C_1/C_0 &= (\sum_{i,j \in I} (D^i \wedge_C D^j))/C_0\\
& =\sum_{i,j \in I} ((D^i \wedge_C D^j)+ C_0)/C_0\\
& \hookrightarrow \sum_{i,j \in I}{^i (C_1/C_0) ^j}\\
& =\bigoplus_{i,j \in I} {^i (C_1/C_0)^j}.\end{align*} This forces the
embedding $((D^i \wedge_C D^j)+C_0)/ C_0
\hookrightarrow {^i(C_1/C_0)^j}$ to be an isomorphism, and hence

\begin{align*}C_1/C_0 &= \sum_{i,j \in I} ((D^i \wedge_C D^j)+ C_0)/C_0\\
& = \bigoplus_{i,j \in I}((D^i \wedge_C D^j)+ C_0)/C_0\\&
\cong \bigoplus_{i,j \in I}(D^i \wedge_C D^j)/ (D^i + D^j).\ \ \ \ \square\end{align*}

\vskip10pt

Theorem 4.1 also permits us to slightly modify the definition of the link quiver of a coalgebra, by adding
multiples of arrows. Of course, in the case for basic coalgebras, it is exactly the original definition.

\vskip 10pt

\begin{defn} ([M2], 1.1) Let $C$ be a coalgebra. The link quiver of $C$ is defined as follows. The vertices are
the isoclasses of simple
subcoalgebras of $C$; and for two simple
subcoalgebras $D^i$ and $D^j$ of $C$, there are
$$l_{ij}:= \frac{1}{n_in_j}{\rm dim}_{K}(D^j \wedge_C D^i)/ (D^i + D^j)$$
arrows from $i$ to $j$, where $n_i$ is the positive integer such that
$(D^i)^* \simeq
M_{n_i}(\Delta_i)$, where $\Delta_i$ is a finite-dimensional
division algebra over $K$.
\end{defn}

\vskip 10pt

\begin{cor} The link quiver of coalgebra $C$ coincides with the Gabriel quiver of $C$.
\end{cor}
\noindent{\bf Proof} \ \ It follows from Theorem 4.1(iv) that

\begin{align*}l_{ij}&= \frac{1}{n_in_j}{\rm dim}_{K}(D^j \wedge_C D^i)/ (D^i + D^j)\\& =
\frac{1}{n_in_j} {\rm dim}_{K} \ {^j(C_1/C_0)^i}\\& =
{\rm dim}_{K} \ e_i.{^j(C_1/C_0)^i}.e_j\\& = t_{ij}. \ \ \ \square
\end{align*}

\vskip10pt

\section{Locally finite coalgebras}

\vskip 10pt

We give a new characterization of locally finite coalgebras, as an
application of Theorems  3.1 and 4.1.

\subsection{} By definition a coalgebra $C$ is said to be locally finite,
provided that the wedge $V \wedge_C W$ is finite-dimensional
whenever $V$ and $W$ are both finite-dimensional. By the
fundamental theorem on coalgebras (i.e., each finite-dimensional
subspace
 of a coalgebra is contained in a finite-dimensional
subcoalgebra), it is clear that a coalgebra $C$ is locally finite
if and only if $D \wedge_C D$ is finite-dimensional for each
finite-dimensional subcoalgebra $D$ of $C$.

Heyneman-Radford showed that a reflexive coalgebra is locally finite
([HR], 3.2.4); Conversely, if $C$ is locally finite with $C_0$
finite-dimensional, then $C$ is reflexive ([HR], 4.2.6).

\vskip10pt

Recall that a subcoalgebra $D$ of $C$ is said to be saturated
provided that $D \wedge_C D=D$.

\vskip5pt

Let $C= C' \oplus C''$ as coalgebras and let $ (M, \rho_l,\rho_r)$
be a $C$-$C$-bicomodule. Set $$N= \{ m \in M \ |\  \rho_l(m) \in
C' \otimes M, \ \ \rho_r(m) \in M \otimes C'\}.$$ Then $N$ is a
$C'$-$C'$-bicomodule.

\begin{lem} With the notation above, ${\rm Cot}_{C'} (N)$ is a saturated
subcoalgebra of ${\rm Cot}_{C} (M)$. \end{lem}

\noindent{\bf Proof} \ \ Set $\widetilde{C}: = {\rm Cot}_{C} (M)$.
By construction of ${\rm Cot}_{C'} (N)$ we have

$${\rm Cot}_{C'} (N)= \bigcup_{n\ge
1}\wedge_{\widetilde{C}}^n
  C'$$
where $\wedge_{\widetilde{C}}^n C'=C'\wedge_{\widetilde{C}}
C'\wedge_{\widetilde{C}}\cdots \wedge_{\widetilde{C}} C'$ ($n$
times). Hence ${\rm Cot}_{C'} (N)$ is saturated in $\widetilde{C}$
(see [HR], 2.1.1). \hfill $\square$

\vskip10pt

The main result of this section is

\vskip10pt

\begin{thm} The Gabriel quiver of a locally finite coalgebra $C$ is locally finite
(i.e., there are only finitely many arrows between arbitrary two
vertices).

Conversely, if the Gabriel quiver of $C$ is locally finite and
$C_0$ is separable, then $C$ is locally finite.
\end{thm}

\noindent{\bf Proof} \quad The necessity follows form Corollary 4.4, since simple coalgebras
are finite-dimensional by the fundamental theorem of coalgebras.

\vskip10pt

Conversely, assume that the Gabriel quiver of $C$ is locally
finite and $C_0$ is separable. In order to prove that $C$ is
locally finite, by Theorem 3.1 it suffices to show that the
cotensor coalgebra ${\rm Cot}_{C_0} (C_1/C_0)$ is locally finite.
This is because a subcoalgebra of a locally finite coalgebra is
again locally finite ([HR], 2.3.2).  In the following, we denote
${\rm Cot}_{C_0} (C_1/C_0)$ by $\widetilde{C}$.\vskip 5pt

Let $D$ be an arbitrary  finite-dimensional subcoalgebra  of
$\widetilde{C}$.  Then the coradical $D_0$ of $D$ is a direct
summand of $C_0$. Set
 \begin{center}
 $M:=\{x \in C_1/C_0 \ |\  \rho_l(x) \in D_0 \otimes C_1/C_0, \ \  \rho_r(x) \in C_1/C_0 \otimes
 D_0\}$.
 \end{center}
Note that $D_0$ is finite-dimensional and $M$ is contained in a
direct sum of finitely many $^i (C_1/C_0)^j$'s. Since the Gabriel
quiver of $C$ is locally finite, it follows that each $^i
(C_1/C_0)^j$ is finite-dimensional, and hence $M$ is
finite-dimensional. By Theorem 3.1 we have
$$D\subseteq {\rm Cot}_{D_0} (D_1/D_0)\subseteq {\rm Cot}_{D_0}
(M).$$ By Lemma 5.1 ${\rm Cot}_{D_0} (M)$ is a saturated
subcoalgbra of $\widetilde{C}$. It follows that
$$D\wedge_{\widetilde{C}}D\subseteq {\rm Cot}_{D_0}
(M)\wedge_{\widetilde{C}}{\rm Cot}_{D_0} (M) = {\rm Cot}_{D_0}
(M).$$ Since $D$ is of finite dimension, we may assume that $D
\subseteq \oplus_{i \leq n} M^{\square i}$ for some $n$. It
follows that $D \wedge_{\widetilde{C}} D$ is contained in
$\oplus_{i \leq 2n} M^{\square i}$ (see Remark 1.1), which is also of finite
dimension. This proves that the cotensor coalgebra ${\rm
Cot}_{C_0} (C_1/C_0)$ is locally finite. \hfill $\square$

\vskip10pt

\section{Quasi-coFrobenius coalgebras}

\vskip10pt

\subsection{} Recall that a coalgebra $C$ is
said to be left quasi-coFrobenius  if there exists an injective
$C^*$-module map from $C$ to a free left $C^*$-module, where the
left $C^*$-module structure on $C$ is given as in 4.2. Similarly,
one has the concept of right quasi-coFrobenius coalgebras. A
coalgebra is quasi-coFrobenius if it is both left
quasi-coFrobenius and right quasi-coFrobenius.  Note that a
coalgebra $C$ is left quasi-coFrobenius if and only if every
injective right $C$-comodule is projective ([DNR], Theorem 3.3.4);
and that if $C$ is left quasi-coFrobenius, then $C^*$ is right
quasi-Frobenius ([DNR], Corollary 3.3.9). Also note that if $C$ is
finite-dimensional, then $C$ is left quasi-coFrobenius if and only
if $C$ is right quasi-coFrobenius, if and only if $C^*$ is
quasi-Frobenius.

\vskip 10pt

We need the following fact, which seems to be well-known.

\vskip10pt

\begin{lem} Let $C$ be a coalgebra. Then $C$ is indecomposable if and only if
the dual algebra $C^*$ is indecomposable.
\end{lem}

\noindent{Proof} \quad Note that here $C$ is not necessarily
finite-dimensional.  The ``if" part is trivial. It suffices to
prove the ``only if " part. If $C^* \simeq A_1 \times A_2$ as
algebras, then ${C^*}^\circ \simeq A_1^\circ \oplus A_2^\circ$ as
coalgebras, where $A^\circ$ denotes the finite dual of  algebra
$A$. Let $\phi: C \longrightarrow {C^*}^*$ be the natural
embedding. Then the image of $\phi$ is contained in ${C^*}^ \circ$
(see e.g. [DNR], Proposition 1.5.12). Identify $C$ with $\phi(C)$.
Then
\begin{center}$C \simeq
({A_1}^\circ \cap C) \oplus ({A_2}^\circ \cap C)$
\end{center} as
coalgebras. Note that  $A_i^\circ \cap C \neq \{0\}$ ($i=1,2$)
(Otherwise, say $A_1^\circ \cap C =\{0\}$, then $C= A_2^\circ \cap
C$, i.e., $C$ is contained in $A_2^{\circ}$, it follows that $A_1$
vanishes on $C$, and hence $A_1 = 0$). This completes the proof.
\hfill $\square$

\vskip 10pt

 The main result of this section is

\vskip 5pt \vskip 5pt

\begin{thm}
Let $C$ be an indecomposable non-simple coalgebra. If $C$ is a
left quasi-coFrobenius, then the Gabriel quiver of $C$ has no
sources.

Thus,  the  Gabriel quiver of a non-simple quasi-coFrobenius coalgebra has no
sources and no sinks.
\end{thm}

\noindent{Proof} \quad  Otherwise, assume that the Gabriel quiver
$Q$ of $C$ has a source $i \in I$. Let $S_i$ be the corresponding
right simple comodule. Then by Lemma 2.3 we have
$soc(E(S_i)/{S_i})\simeq \oplus_{j \in I} \frac{t_{ji}}{d_j}S_j$.
Since $i$ is a source, it follows that  $t_{ji}=0$ for every  $j
\in I$, and hence $E(S_i)=S_i$.

\vskip10pt

For any $j \neq i$, we have
$${\rm Hom}_C(E(S_i), E(S_j))= {\rm Hom}_C(S_i, E(S_j))= {\rm
Hom}_C(S_i, S_j)=0,$$ here we have used $soc(E(S_j))=S_j$ and the
Schur lemma.

On the other hand, since $C$ is left quasi-coFrobenius, it follows
that $E(S_i) = S_i$ is projective as right $C$-comodule. Thus
$${\rm Hom}_C(E(S_j), E(S_i))=0$$ for each $
 j \neq i  \in I $ (otherwise, let $f: E(S_j) \longrightarrow
E(S_i)=S_i$ be a nonzero $C$-comodule map. Then $f$ is surjective
. Thus $E(S_j) \simeq S_i \oplus {\rm Ker}(f)$ by the projectivity
of $S_i$, a contradiction).

Note that $C$ itself is a right $C$-comudule via $\Delta_C$, and
that there is an algebra isomorphism ${\rm End}_C(C) \simeq C^*$,
sending $f$ to $ \varepsilon_C \circ f$ (see e.g.  [DNR],
Proposition 3.1.8). Since $C \simeq \oplus_{j \in I} n_jE(S_j)$ as
right $C$-comodule, it follows that
$$C^* \simeq {\rm End}_C (C)
\simeq {\rm End}_C(n_iE(S_i)) \oplus {\rm End}_C(\oplus_{j \neq i}
(n_j E(S_j))).$$ While $C^*$ is indecomposable by Lemma 6.1, we
then obtain a desired contradiction.  \hfill $\square$

 \vskip 20pt
\bibliography{}

\end{document}